\documentclass[11pt, epsfig]{article}
\usepackage{epsfig, amsmath, amssymb, amsthm, times}
\usepackage{graphicx}

\usepackage[mathcal]{eucal}
\usepackage{bbm}

\newtheorem {thm}{Theorem}[section]
\newtheorem {lem}[thm]{Lemma}

\theoremstyle{defintion}
\newtheorem {df}[thm]{Definition}

\theoremstyle{remark}
\newtheorem{rem}[thm]{Remark}

\theoremstyle{example}
\newtheorem{ex}[thm]{Example}

\theoremstyle{assumption}

\def\pf{{\it Proof.\;}}

\def\E{{\mathbb E~}}
\def\P{{\mathbb P~}}
\def\R{{\mathbb R}}
\def\N{{\mathbb N}}

\def\lbl{\label}
\def\be{\begin{equation}}
\def\ee{\end{equation}}
\def\p{\partial}
\def\qed{\square}

\def\deg{\operatorname{deg}}

\def\one{\mathbf{1}}

\title{The semilinear heat equation on sparse random graphs}
\author{Dmitry Kaliuzhnyi-Verbovetskyi and Georgi S. Medvedev
\thanks{
Department of Mathematics, Drexel University, 3141 Chestnut Street,
Philadelphia, PA 19104;
{\tt dmitryk@math.drexel.edu},
{\tt medvedev@drexel.edu}
}
}

\begin{document}
\maketitle
\begin{abstract}
Using the theory of $L^p$--graphons \cite{BCCZ-I-14, BCCZ-II-14}, we derive and 
rigorously justify the continuum limit for systems of differential equations on sparse
random graphs. Specifically, we show that the solutions of the initial value problems 
for the discrete models can be approximated by those of an appropriate nonlocal 
diffusion equation. Our results apply to a range of spatially extended dynamical models
of different physical, biological, social, and economic networks.
Importantly, our assumptions cover network topologies  featured in many important
real-world networks. In particular, we derive  the continuum limit for  coupled
dynamical systems on power law graphs. The latter is the main
motivation for this work. 
\end{abstract}

\section{Introduction}
Reaction-diffusion equations describe the change of concentration of
chemical, biological, or other species as a result of local reaction
and spatial diffusion:
\be\lbl{RD}
{\p\over\p t} u =\Delta u + f(u).
\ee
Here, $u:Q\times \R^+$ is an unknown function, whose interpretation depends on the
model at hand, defined on spatial domain $Q\subset\R^n$ and evolving in time.
Reaction-difussion systems have been successfully used to study pattern formation and
propagation phenomena in such diverse areas of science as ecology,
molecular biology, morphogenesis, neuroscience, and material science,
to name a few \cite{CroHoh93}.

In many models of collective behavior of discrete agents, one
is led to replace the spatial domain $Q$ by a graph and the Laplace operator
$\Delta$ by the graph Laplacian \cite{Biggs}. Specifically, let
$\Gamma_n=\langle V(\Gamma_n), E(\Gamma_n)\rangle$ denote a graph
on $n$ nodes. Here, $V(\Gamma_n)$ and  $E(\Gamma_n)$ stand for the
sets of nodes and edges respectively. Without loss of generality,
let $V(\Gamma_n)=\{1,2,\dots, n\}=:[n]$ and
consider the following nonlinear evolution equation on $\Gamma_n$:
\be\lbl{nlin}
\dot u_{ni} = {1\over \deg_{\Gamma_n}(i)}
\sum_{j: \{i,j\}\in E(\Gamma_n)}  D(u_{nj}-u_{ni}) + f(u_{ni}), \; i\in [n],
\ee
where
$D$ and $f$ are Lipschitz continuous functions and $\deg_{\Gamma_n}(i)$ stands for the
degree of node $i\in [n].$ The sum on the right-hand side of
\eqref{nlin} models  nonlinear diffusion on $\Gamma_n$.
Discrete diffusion operators of this form have been used in various models of collective
behavior. For instance, with $D(u)=\sin u$ it appears in the Kuramoto model
\cite{Kur84} and in the power network models \cite{DorBul12}, with $D(u)=\phi(|u|) u$
for an appropriate function $\phi$, it is used in models of flocking
\cite{CucSma07} 
and  opinion dynamics \cite{MotTad14}, and with $D(u)=u$ - in consensus protocols
\cite{Med12}. In the latter case,  \eqref{nlin} becomes a semilinear heat equation
on $\Gamma_n$:
\be\lbl{semi}
\dot u_{ni} = {1\over \deg_{\Gamma_n}(i)}
\sum_{j: \{j,i\}\in E(\Gamma_n)}  (u_{nj}-u_{ni}) + f(u_{ni}), \; i\in [n].
\ee

Understanding the dynamics of coupled systems \eqref{nlin} and
\eqref{semi} on graphs modeling connectivity in real-life systems like
neuronal networks, power grids, or the Internet, can be quite challenging.
Recently, new powerful techniques for describing and analyzing the
structure of large  graphs, based on the appropriate notions of convergence, 
have been developed in the graph theory \cite{LovGraphLim12}.
Many nontrivial graph sequences that are of interest in applications, such
as Erd\H{o}s-R{\' e}nyi, small-world, and preferential attachment graphs, 
have relatively simple limits, expressed by symmetric measurable functions 
on a unit square, called graphons \cite{LovGraphLim12}.  These graph limits 
can be used for developing continuum models approximating 
the dynamics of \eqref{nlin} for large $n$:
\be\lbl{climit}
{\p\over \p t} u(x,t)= \int_I W(x,y) D( u(y,t)-u(x,t)) dy + f(u(x,t)),
\ee
where $W$ is the graphon describing the limiting behavior of $\{\Gamma_n\}.$

In \cite{Med14a, Med14b}, the continuum limit \eqref{climit} was derived and rigorously
justified for coupled dynamical systems on convergent families of dense 
graphs\footnote{If $|E|=O(|V|^2),$ $\Gamma=\langle V, E\rangle$ is called dense,
otherwise it is called sparse.}. The analysis in \cite{Med14a, Med14b}
covers systems on many interesting graphs including small-world and Erd\H{o}s-R{\' e}nyi graphs. 
However, many real-world networks feature sparse connectivity.
Thus, in this paper, we focus on coupled systems on convergent families of 
sparse graphs. 

Our work is inspired by the recent progress made by
Borgs, Chayes, Cohn, and Zhao, who extended the theory of graph limits originally developed
for dense graphs to  sparse graphs of unbounded degree \cite{BCCZ-I-14,BCCZ-II-14}.
The new theory covers many interesting examples of graphs. Notably, it
applies to graphs with power law degree distribution,  which was identified in many different systems \cite{BarAlb99}. 
A distinctive feature of the convergence theory for sparse graphs is that the graphons are 
no longer  bounded, as in the dense case, but in general are functions from $L^p(I^2)$, $p>1$. 
This leads to continuum model \eqref{climit} with $W\in L^p(I^2)$. The analysis of \eqref{climit}
with an $L^p$ kernel presents new challenges that were not present in the $L^\infty$-case, analyzed
in \cite{Med14a, Med14b}. Overcoming these problems is the goal of this paper.  

In the next section, we adapt the notion of W-random sparse graph from \cite{BCCZ-I-14} to 
define a  sequence of random graphs $\Gamma_n=G(W,\rho_n,X_n)$ with edge density 
$\rho_n\to 0$ as $n\to\infty$ and with the graph limit $W\in L^p(I^2), p>1.$
This random graph model will be used throughout this paper. It covers power law graphs, 
our main motivating example, as well as sparse stochastic block and sparse 
Erd\H{o}s-R{\' e}nyi graphs (cf.~Examples~\ref{Walpha} and \ref{ex.sparse}) 
among many other sparse graphs. In \S\ref{sec.graphs}, we compute the expected degree
and edge density of $\Gamma_n=G(W,\rho_n,X_n)$. We then formulate the dynamical model on 
$\Gamma_n$ and formally derive the continuum limit \eqref{climit}. The derivation
includes two steps. First, we average the right hand side of the coupled model (which depends
on the random realization of $\Gamma_n$) to obtain a deterministic equation.
We then send $n\to \infty$ in the averaged system to derive the continuum limit.
This derivation is done for the semilinear heat equation \eqref{semi}, which will be studied 
in the main part of the paper. However, the same derivation easily translates to the nonlinear 
equation \eqref{nlin}, which results in \eqref{climit}. In Section~\ref{sec.ivp}, we establish well-posedness
of the IVP for \eqref{climit} and derive certain a priori estimates for solutions of the initial value problems (IVPs). 
For technical
reasons, we restrict to studying the semilinear heat equation for
the remainder of this paper. In the last section, we comment on how
this analysis extends to cover to certain nonlinear models arising in
applications. In particular, we discuss the Kuramoto
model on power law graphs.

The main result of this work is formulated in \S\ref{sec.main}. 
Under the appropriate assumptions on the graphon $W\in L^2(I^2)$ and the nonlinearity
$f$, we prove that the solutions of the IVP for the semilinear heat equation \eqref{semi}
on $\Gamma_n$ converge in $L^2(0,T; L^2(I))$ (for any  $T>0$) almost surely (a.s.) to the solution of the continuum limit 
\eqref{climit} subject to the appropriate initial condition as $n\to\infty$. This is the subject of 
Theorem~\ref{thm.main}, which is proved in Sections~\ref{sec.averaging} and \ref{sec.climit}.
In the former section, the justification for the averaging is provided. In the latter, we show that the solutions
of the IVP for the averaged equation on $\Gamma_n$ converge to those for the continuum limit 
as $n\to\infty$. To this end, we show that the solutions of the averaged equation can be approximated
by the solutions of certain Galerkin problems, which, in turn,  converge to the solution of the continuum
limit. The final section discusses extensions of our work to certain nonlinear models that are important in
applications.

\section{The model}\lbl{sec.model}
\setcounter{equation}{0}

\subsection{The random graph model} \lbl{sec.graphs}
We start with the description of the sparse random  graphs that will be used in
this paper. Our random graph model is motivated by the construction of
sparse $W$-random
graphs in \cite{BCCZ-I-14,BCCZ-II-14}. Specifically, let $W$ be a symmetric
nonnegative function on a unit square $I^2$, $X_n$ be a discretization of
$I$
\be \lbl{Xn}
X_n = \{ x_{n0}, x_{n1},x_{n2},\dots,x_{nn}\}, \; x_{ni}=i/n,\; i=0,1,\dots, n,
\ee
and $\{\rho_n\}$ be a sequence of positive numbers such that
$\rho_n\to 0$ and $n\rho_n\to\infty$ as $n\to\infty$.

$\Gamma_n=G(W,\rho_n,X_n)$ stands for a random graph with the node
set $V(\Gamma_n)=[n]$ and the edge set $E(\Gamma_n)$ such that
the probability that $\{i, j\}$ forms an edge is
\be\lbl{Pedge}
\P(\{i,j\} \in E(\Gamma_n))=\rho_n \bar W_n(x_{ni},x_{nj}), \; i,j\in [n],\footnotemark
\ee
\footnotetext{To keep notation simple, we allow for loops, i.e., edges connecting a node
to itself, in our random graph model. Excluding loops would not lead
to any changes in the analysis.}
where
\be
\lbl{Wn}
\bar W_n(x,y)= \rho_n^{-1}\wedge W(x,y).\footnotemark
\ee
\footnotetext{Throughout this paper,
  we use $a\wedge b$ and $a\vee b$ to denote $\min\{a,b\}$ and
  $\max\{a,b\}$ respectively.}
The decision whether a given pair of nodes is included in the edge set
is made independently from other pairs. In other words,
$G(W,\rho_n,X_n)$ 
is a product probability
space
\be\lbl{pspace}
(\Omega_n= \{0,1\}^{n(n+1)/2}, 2^{\Omega_n},\P).
\ee
By $\Gamma_n(\omega), \omega\in \Omega_n$, we will denote
a random graph drawn from the probability distribution $G(W,\rho_n,X_n)$.

Throughout this paper, we use Bernoulli random variables
\be\lbl{xi}
\xi_{nij}(\omega)=\one_{\{i,j\}\in E(\Gamma_n)}(\omega), \; i,j\in [n].
\ee
to describe the edge set of $\Gamma_n$.
Random variable $\xi_{nij}$  takes value $1$ if $\{i,j\}$ forms an edge and $0$ otherwise. In particular,
\be\lbl{Exi}
\E\xi_{nij}= \P(\{i,j\}\in E(\Gamma_n))=\rho_n\bar W_n(x_{ni},x_{nj}), \; \{i,j\}\in [n],
\ee
and the expected degree of node $i$ of $\Gamma_n$
\be\lbl{exp-deg}
\E\operatorname{deg}_{\Gamma_n}(i)=
\E\left(\sum_{j=1}^n \xi_{nij}\right)=\rho_n\sum_{j=1}^n \bar W_n(x_{ni},x_{nj}).
\ee

Next, we formulate our assumptions on the graphon $W$.
\begin{description}
\item[W-1)] $W\in L^2(I^2)$ is a nonnegative symmetric on the unit square
$I^2$ that   is continuous on its interior.
\item[W-2)] $\int_{I^2} W(x,y) dxdy>0$ and
\be\lbl{W-2}
 n^{-2}\sum_{i,j=1}^n \bar W_n(x_{ni},x_{nj}) =\int_{I^2} W(x,y) dxdy +o(1).
\ee
\item[W-3)] For every $x\in (0,1],$ $W(x,\cdot)\in L^1(I),$ and
\be\lbl{int-below}
\inf_{x\in (0,1)} \int_I W(x,z)dz =:\nu >0.
\ee
Moreover,
\be\lbl{W-3} 
 n^{-1}\sum_{j=1}^n \bar W_n(x,x_j) =\int_{I} W(x,y) dy\left(1+\delta_n(x)\right),
\ee
where $\delta_n(x)\to 0$ as $n\to\infty$ uniformly in $x\in (0,1)$.
\end{description}

Conditions in \textbf{W-2)} and \textbf{W-3)}
guarantee that the expected edge density and expected degrees of nodes of $\Gamma_n$
for $n\gg 1$ are well-defined and are well-approximated by the corresponding integrals of $W$.
The above assumptions on graphon $W$ are dictated by the random graph model and are
practically minimal.

 We will now introduce two more technical assumptions that are needed
 for the proof of our main result:
 \noindent
 \textbf{W-4)} $W\in L^4(I^2)$ and 
 \be\lbl{W-4}
 \limsup_{n} n^{-2} \sum_{i,j=1}^n \bar W_n(x,x_j)^2< \infty.
 \ee
 Assumptions in \textbf{W-4)} will not be used until \S\ref{sec.approximate}.

The main example motivating our random graph model is the following configuration model
of a power law graph.
\begin{ex}\lbl{ex.pl} Let $0<\alpha<\gamma<1$ and consider
$G(W,\rho_n,X_n)$, where $\rho_n=n^{-\gamma}$ and
\be \lbl{Walpha}
W(x,y)= (1-\alpha)^2 (xy)^{-\alpha}.
\ee
\end{ex}

\begin{lem}\lbl{lem.pl} $\Gamma_n=G(W,\rho_n,X_n)$ of Example~\ref{ex.pl} is a sparse
graph with a power law expected degree distribution. In particular, we have
\begin{description}
\item[A)] The expected degree of node $i\in [n]$ of $\Gamma_n$ is
\be\lbl{Edeg}
\E\deg_{\Gamma_n}(i)=(1-\alpha)n^{1+\alpha-\gamma} i^{-\alpha}(1+o(1)).
\ee
\item[B)] The expected edge density of $\Gamma_n$ is $n^{-\gamma}(1+o(1)).$
\end{description}
\end{lem}
\pf
By \eqref{Pedge},
\be\lbl{Edeg-1}
\E\deg_{\Gamma_n}(i)=\sum_{j=1}^n \P(\{i,j\}\in E(\Gamma_n))=
\rho_n n \sum_{j=1}^n \bar W_n(x_{ni},x_{nj}) n^{-1}.
\ee
Plugging  $\rho_n=n^{-\gamma}$ and  (\ref{Wn}) in \eqref{Edeg-1}, we have as
\be\lbl{Edeg-2}
\E\deg_{\Gamma_n}(i)=(1-\alpha) n^{1+\alpha-\gamma} i^{-\alpha}
\sum_{j=1}^n \left[\delta_n \wedge W^{(1)} (x_{nj}) \right] n^{-1},
\ee
where $W^{(1)}(x)=(1-\alpha) x^{-\alpha}$ and $\delta_n=i^\alpha n^{\gamma-\alpha} (1-\alpha)^{-1}.$

Denote
\be\lbl{df-Ini}
I_{ni}:=(x_{n(i-1)}, x_{ni}], \quad i\in [n].
\ee
Next, let
\be\lbl{def-step-W1}
W_n^{(1)}(x)=\sum_{j=1}^n (\delta_n\wedge W^{(1)}(x_{nj})) \one_{I_{nj}} (x)
\ee
and note that
\be\lbl{integral-of-step-W}
n^{-1}\sum_{j=1}^n \delta_n \wedge W_n^{(1)} (x_{nj}) =\int_I W_n^{(1)} dx.
\ee
Furthermore, $W_n^{(1)}\le W^{(1)}$ and $W_n^{(1)}\to W^{(1)}$
pointwise on $(0,1]$ as $n\to\infty$.
By the Dominated Convergence Theorem \cite{Bogachev-MT},
\be\lbl{take-the-limit}
\lim_{n\to\infty} \int_I W_n^{(1)} dx=\int_I W^{(1)} dx=1.
\ee
The combination of \eqref{Edeg-2}, \eqref{integral-of-step-W}, and \eqref{take-the-limit}
yields \eqref{Edeg}. This shows \textbf{A)}.

A similar argument is used to estimate the expected number of edges in $\Gamma_n$
\begin{eqnarray*}
\E |E(\Gamma_n)| &=& {1\over 2} \sum_{i=1}^n\sum_{j=1}^n \rho_n \bar W_n(x_{ni},x_{nj})
= {1\over 2} n^{2-\gamma}\sum_{i=1}^n\sum_{j=1}^n  \bar W_n(x_{ni},x_{nj}) n^{-2}.
\end{eqnarray*}
Define
\be\lbl{2d-step}
W_n(x,y)=\sum_{i,j=1}^n (\rho_n^{-1}\wedge W(x_{ni},x_{nj}))
 \one_{I_{ni}\times I_{nj}}(x,y).
\ee Then \be\lbl{relateto-volume} \E|
E(\Gamma_n)|={n^{2-\gamma}\over 2} \int_{I^2} W_n dxdy. \ee
By construction, $W_n\le W$ and $W_n\to W$  as $n\to\infty$ on
$(0,1]\times (0,1].$ By the Dominated Convergence Theorem,
\be\lbl{DC-again} \lim_{n\to\infty}\int_{I^2} W_n dxdy =\int_{I^2}
W dxdy=1. \ee Equations \eqref{relateto-volume} and
\eqref{DC-again} imply \be\lbl{exp-volume} \E|
E(\Gamma_n)|={n^{2-\gamma}\over 2}(1+o(1)). \ee By dividing both
sides of \eqref{exp-volume} by $n(n+1)/2,$ the total  number of
possible edges, we obtain the statement in \textbf{B)}. $\qed$

\begin{rem}
The power law graphs defined  above are sparse, because the
expected edge density is $O(n^{-\gamma})$, $\gamma>0.$ On the
other hand, the expected number of edges grows superlinearly
$n^{2-\gamma}, $ because $\gamma<1$. To preserve these features,
in the general random graph model  $G(W,\rho_n, X_n), n\in\N,$
above it is assumed that $\rho_n\to 0$ and $n\rho_n\to\infty$ as
$n\to\infty.$
\end{rem}

We conclude the discussion of the graph model with two more examples of sparse graphs
covered by our assumptions. Both examples are taken from
\cite{BCCZ-II-14}.

\begin{ex}\lbl{ex.sparse} Consider $\Gamma_n=G(W,n^{-\beta},X_n) , \; \beta\in (0,1)$
 for the following choices of $W$.
\begin{description}
\item[1)]
Let $W\equiv 1.$ Then $\Gamma_n$ is a generalization of an
Erd\H{o}s-R\'{e}nyi random graph $G_{n,p}$ with $p=n^{-\beta}$. Note that the edge density
in this case is $n^{-\beta}$. For the classical Erd\H{o}s-R\'{e}nyi graph $G_{n,p}$
with constant $p\in (0,1)$, the edge density is equal to $p$.
The latter graph is dense, whereas the former is sparse.
\item[2)]
Let $W(x,y)=b_{ij}\ge 0,$ $(x,y)\in V_i\times V_j, (i,j)\in [n]^2,$
where $\sum_{i,j=1}^k b_{ij}>0$ and
$(V_1, V_2, \dots, V_k)$ is a partition of $I$ into $k$ disjoint intervals.
In this case, $\Gamma_n$ is a sparse stochastic block graph with
edge density $n^{-\beta}$.
\end{description}
\end{ex}

\subsection{The dynamical model}\lbl{sec.dynamics}
Having defined the structure of the network, we next turn to its  dynamics.
Let $\Gamma_n=\Gamma_n(\omega), \; \omega\in\Omega_n$ (cf. \eqref{pspace})
be a random graph taken from the probability
distribution $G(W,\rho_n, X_n)$
and consider the following system of differential equations
\be\lbl{heat}
\dot u_{ni}={1\over d_{ni}}
\sum_{j=1}^n \xi_{nij}(\omega)
\left(u_{nj}-u_{ni}\right) + f(u_{ni}), \; i\in [n],
\ee
where  $u_n(t)=(u_{n1}(t),u_{n2}(t),\dots, u_{nn}(t)),$
$\xi_{nij}, \; i,j\in [n],$ are Bernoulli random variables defined in \eqref{xi},
$
d_{ni}=\E\deg_{\Gamma_n}(i)
$
(cf. \eqref{exp-deg}), and
$f:\R\to\R$ is
a Lipschitz
continuous function
\be\lbl{Lip-f}
\exists  L_f>0: \quad  |f(x)-f(y)|\le L_f |x-y|,\; \forall x,y\in\R.
\ee

The sum on the right hand side of \eqref{heat} defines a discrete diffusion
operator.
For simplicity, we scale the sum
on the right-hand side of \eqref{heat} by the expected degree rather than by the actual degree.
Thus, \eqref{heat} is an evolution equation on  a random graph $\Gamma_n$.
Specifically, it is  a semilinear heat equation on  $\Gamma_n,$
since the sum on the right-hand side of \eqref{heat} is a discrete graph Laplacian.

We are interested in describing the dynamics of \eqref{heat}  for $n\gg 1$.
However, the right-hand side of (\ref{heat}) depends on the random graph
$\Gamma_n(\omega)$, i.e., on the random event $\omega\in\Omega_n$:
$$
F_{ni}(u_n,\omega)={1\over d_{ni}}\sum_{j=1}^n
\xi_{nij}(\omega) \left(u_{nj}-u_{ni}\right) +f(u_{ni}),\quad i\in [n].
$$
As the first step in the analysis of \eqref{heat}, we approximate it by the deterministic
problem
by averaging the right-hand side of \eqref{heat} over all realizations of $\Gamma_n$:
\be\lbl{aKM}
 \dot v_{ni}(t)=\bar F_{ni}(v_n), \; v_n(t)=(v_{n1}(t),v_{n2}(t),\dots, v_{nn}(t)),
\ee
where
\begin{eqnarray} \nonumber
\bar F_{ni} (v)&=& \E F_{ni} (v_n,\cdot)=
{\rho_n\over d_{ni}} \sum_{j=1}^n \bar W_n(x_{ni}, x_{nj}) (v_{nj}-v_{ni}) +f(v_{ni}), \\
\lbl{average}
&=& {1\over n} \sum_{j=1}^n V_{nij}
(v_{nj}-v_{ni}) +f(v_{ni}), \quad i\in [n],
\end{eqnarray}
where
\be\lbl{def-Vnij}
V_{nij}= {\bar W_n(x_{ni}, x_{nj})\over n^{-1} \sum_{k=1}^n \bar
  W_n(x_{ni}, x_{nk})}.
\ee
Next, we take the limit in the averaged equation \eqref{average} as $n\to\infty.$
To this end, we represent the solution of \eqref{average} as a step function
\be\lbl{v-step}
v_n(x,t)=\sum_{i=1}^n v_{ni}(t) \one_{I_{ni}} (x),\;  
\ee
and rewrite \eqref{average} as
\be\lbl{average-step}
{\p\over \p t} v_n(x,t)=\int_I V_n(x,y)(v_n(y,t)-v_n(x,t)) dy + f(v_n(x,t)),
\ee
where
\be\lbl{W-step}
V_n(x,y)= \sum_{i,j=1}^n V_{nij} \one_{I_{ni}\times I_{nj}}.
\ee
Assuming that $v_n(x,t)\to u(x,t)$ in the appropriate sense, and using
the integrability assumptions \textbf{W-2)} and \textbf{W-3)},
in the limit as $n\to\infty$ we formally obtain the following continuum limit of (\ref{aKM})
\be\lbl{cKM}
{\p \over \p t} u(x,t)=\int_I U(x,y) \left( u(y,t)-u(x,t)\right)dy + f(u(x,t)),
\ee
\be\lbl{kernel}
U(x,y)=W(x,y)\left(\int_I W(x,z)dz\right)^{-1}.
\ee
Note that $U\in L^2(I^2)$ (cf.~\textbf{W-1)} and \eqref{int-below}). 
%
%

\begin{ex}\lbl{ex.continuum-pl}
For the power law graphs defined in Example~\ref{ex.pl} with square integrable graphons,
the continuum limit takes the following form
$$
{\p\over \p t} u (x,t)=\int_I y^{-\alpha} \left( u(y,t)-u(x,t)\right)dy + f(u(x,t)),\; 0<\alpha<1/2.
$$
\end{ex}

The goal of this paper is to describe the relation between the
solutions of the IVPs for the discrete model \eqref{heat} on
sparse graph $\Gamma_n, n\gg 1,$  the averaged model \eqref{aKM}, and
the continuum limit \eqref{cKM}.

\subsection{The main result}\lbl{sec.main}

Let $g\in L^\infty(I)$ and consider the IVP for \eqref{cKM} subject to the 
initial condition
\be\lbl{ic-clim}
u(x,0)=g(x), \quad x\in I.
\ee
Likewise, we supply the discrete problem \eqref{heat} with the initial condition
\be\lbl{ic-heat}
u_{ni}(0)=n\int_{I_{ni}} g(x)dx,\quad i\in [n].
\ee
To compare solutions of the IVPs for the discrete and continuous models, we define the 
step function 
\be\lbl{u-step}
u_n(x,t)=\sum_{i=1}^n u_{ni}(t) \one_{I_{ni}}(x).
\ee
The main result of this paper concerns the $L^2$-proximity of $u_n(\cdot,t)$ and $u(\cdot,t)$
on finite time intervals for large $n$.

\begin{thm}\lbl{thm.main}
Let $\Gamma_n=G(W,\rho_n, X_n)$ be a sequence of random graphs, where $W$ satisfies
conditions
\newline \textbf{W-1)}-\textbf{W-4)}, $X_n, n\in \N,$ are defined in \eqref{Xn} and the positive 
sequence $\{\rho_n\}$ is such that $\rho_n\to 0$ and $n\rho_n\to\infty$ as $n\to\infty$.
Suppose $f:\R\to\R$ is Lipschitz continuous function (cf.~\eqref{Lip-f}), $g\in L^\infty(I)$, and
$T>0$ is arbitrary.

Then with probability $1,$ for solutions of the IVPs \eqref{cKM}, \eqref{ic-clim} and \eqref{heat},
\eqref{ic-heat}, we have
$$
\int_0^T \|u_n(\cdot, t)-u(\cdot, t)\|^2_{L^2(I)} dt \to 0, \quad \mbox{as}\; n\to\infty.
$$
\end{thm}

\section{The IVP for the nonlocal equation}\lbl{sec.ivp}
\setcounter{equation}{0}
\subsection{Existence and uniqueness of solutions}\lbl{sec.exist}

In this section, we show that the IVP for \eqref{cKM}, \eqref{kernel} has a unique solution.
The contraction mapping principle used below applies to the nonlinear
heat equation
\be\lbl{cKM+}
{\p \over \p t} u(x,t)=\int_I U(x,y) D\left( u(y,t)-u(x,t)\right)dy + f(u(x,t)),
\ee
where $D\colon\R\to\R$ is  Lipschitz continuous:
\be\lbl{Lip}
\left|D(u)-D(v)\right|\le L_D |u-v| \;\forall u,v\in \R.
\ee
Below, we study the well-posedness of the IVP for \eqref{cKM+}. The results of this section
will obviously hold for \eqref{cKM} as well.

With the definition \eqref{kernel} in mind,
in this section, we assume that $U\in L^p(I^2), p\ge 2,$ is a
nonnegative function, satisfying
\be\lbl{norm}
\int_I U(x,y)dy=1.
\ee
We interpret the solution of the IVP for (\ref{cKM+}),
$u(x,t)$, as a vector-valued map $\mathbf{u}:\R\to L^q(I)$, i.e.,
$[\mathbf{u}(t)](x)=u(x,t).$

\begin{thm}\lbl{thm.wellposed}
Suppose $U\in L^p(I^2), \, p\ge 2,$ is a nonnegative function
satisfying \eqref{norm} and functions $f$ and $D$ satisfy
\eqref{Lip-f} and \eqref{Lip} respectively. Then the IVP 
for \eqref{cKM+} with initial data $\mathbf{u}(0)=g\in
L^q(I),\; q=p/(p-1)$ has a unique solution 
$\mathbf{u}\in C^1(\R;L^q(I)),$ which depends continuously on 
$g$.
%
\end{thm}
\pf\; Denote 
\be\lbl{tau} \tau= (2L \left(\|U\|_{L^p (I^2)}
+2\right) )^{-1}, \ee where $L=L_f\vee L_D$ is the largest of the
two Lipschitz constants of $D$ and $f$ (cf.~\eqref{Lip-f},
(\ref{Lip})). Denote $\mathcal{M}=C(0,\tau; L^q (I))$ and define
$K\colon\mathcal{M}\rightarrow \mathcal{M}$ as follows:
\be\lbl{df-K} [K\mathbf{u}](t)=g + \int_0^t \left( \int_I
U(\cdot,y) D\left( u(y,s)-u(\cdot,s)\right)dy
+f(u(\cdot,s))\right)ds. 
\ee (The correctness of this definition
will be shown later.) We rewrite the IVP for \eqref{cKM+} as a
fixed point equation for the mapping $K$, 
\be\lbl{fix}
\mathbf{u}=K\mathbf{u}, \ee and show that $K$ is a contraction on
$\mathcal{M}$.

The following inequalities hold for any $u\in L^q(I)$ and $W\in
L^p(I^2), \; p>1,\; q=p/(p-1),$ 
\be\lbl{byHolder}
\|u\|_{L^q(I)}\le \|u\|_{L^{p\vee q}(I)},\; \|W\|_{L^p(I^2)}\le
\|W\|_{L^{p\vee q} (I^2)}. \ee 
They follow from the H\"{o}lder
inequality applied to functions defined on the unit interval $I$
and the unit square $I^2$ respectively. In particular, for $q\le
2\le p,$ we have \be\lbl{as-follows-by-H} \|u\|_{L^q(I)}\le
\|u\|_{L^p(I)},\; \|U\|_{L^q(I^2)}\le \|U\|_{L^p(I^2)}. \ee 
For
any $\mathbf{u}, \mathbf{v}\in \mathcal{M},$ we have
\begin{multline}
\left\| K\mathbf{u} - K\mathbf{v} \right\|_{\mathcal{M}} =
\max_{t\in [0,\tau]} \left\| K\mathbf{u}(t) - K\mathbf{v}(t)
\right\|_{L^q (I)} \\
\le \max_{t\in [0,\tau]} \left\| \int_0^t \left(\int_I U(\cdot,y)
 \left|  D\left( u(y,s)-u(\cdot,s)\right)- D\left(
v(y,s)-v(\cdot,s)\right)
\right| dy   +L|u(\cdot,s)-v(\cdot,s)|\right)  ds \right\|_{L^q (I)}\\
\le L\max_{t\in [0,\tau]} \left\| \int_0^t \left(\int_I U(\cdot,y)
\left|  u(y,s)-u(\cdot,s)- v(y,s)+v(\cdot,s)
\right| dy   +|u(\cdot,s)-v(\cdot,s)|\right)  ds \right\|_{L^q (I)}\\
 \le \tau L \max_{t\in [0,\tau]} \left\{ \left\| \int_I
U(\cdot,y) \left|u(y,t)-v(y,t)\right| dy\right\|_{L^q (I)}+
\left\| \int_I U(\cdot,y) \left|u(\cdot,t)-v(\cdot,t)\right| dy\right\|_{L^q (I)} \right. \\
\lbl{triangle} + \left. \left\| u(\cdot,t)-v(\cdot,
t)\right\|_{L^q (I)} \right\}\\
=\tau L \max_{t\in [0,\tau]} \left\{ \left\| \int_I U(\cdot,y)
\left|u(y,t)-v(y,t)\right| dy\right\|_{L^q (I)}+ 2 \left\|
u(\cdot,t)-v(\cdot, t)\right\|_{L^q (I)} \right\},
\end{multline}
where we used Lipschitz continuity of $D$ and $f$, and
\eqref{norm}. Using the H\"{o}lder inequality and the second
inequality in \eqref{as-follows-by-H}, we have
\begin{eqnarray}\nonumber
\left\| \int_I U(\cdot,y) \left|u(y,t)-v(y,t)\right|
dy\right\|_{L^q (I)} &\le& \left\|
\left\|U(x,\cdot)\right\|_{L^{p} (I)}
\left\| u(\cdot, t)- v(\cdot, t)\right\|_{L^{q} (I)} \right\|_{L^q (I_x)}   \\
\lbl{first-term}
&\le& \|U\|_{L^{p} (I^2)} \left\|\mathbf{u}(t)-\mathbf{v}(t)\right\|_{L^q (I)},
\end{eqnarray}
where $I_x=[0,1]$ refers to the domain of a function of $x$.

The combination of (\ref{triangle}) and (\ref{first-term}) yields
\be\lbl{contraction} \left\| K\mathbf{u} - K\mathbf{v}
\right\|_{\mathcal{M}}\le L\tau \left(\|U\|_{L^{p\vee q} (I^2)}
+2\right) \left\|\mathbf{u}-\mathbf{v}\right\|_{\mathcal{M}}. \ee
Thus, using (\ref{tau}), we have
$$
\left\| K\mathbf{u} - K\mathbf{v} \right\|_{\mathcal{M}}\le
{1\over 2} \left\|\mathbf{u}-\mathbf{v}\right\|_{\mathcal{M}}.
$$
It follows that $K$ is a correctly defined contraction on
$\mathcal{M}$.

Next, we show $K(\mathcal{M})\subset \mathcal{M}$. To this end, for $\mathbf{z}\equiv 0$
on $I\times [0,\tau]$, we have
\begin{eqnarray} \nonumber
\|K\mathbf{u}\|_\mathcal{M}&\le&  \|K\mathbf{u}-K\mathbf{z}\|_\mathcal{M}
+\|K\mathbf{z}\|_\mathcal{M}\\
\lbl{Ku}
&\le& {1\over 2} \|\mathbf{u}\|_\mathcal{M} +\|K\mathbf{z}\|_\mathcal{M}.
\end{eqnarray}
Further,
$$
[K\mathbf{z}](t)=g+t\left(D(0)+ f(0)\right),
$$
so that $Kz\in\mathcal{M}$, and then \eqref{Ku} implies that
$K\mathbf{u}\in\mathcal{M}$.

From \eqref{contraction}, by the Banach contraction mapping principle, there exists a unique
solution of the IVP for (\ref{cKM+})
$\bar{\mathbf{u}}\in \mathcal{M}\subset C(0,\tau; L^q(I))$.
Using $\bar{\mathbf{u}}(\tau)$ as the initial condition, the local solution
can be extended  to $[0, 2\tau]$, and by repeating this argument
to $[0, T]$ for any $T>0$. In a similar fashion, we can prove the existence and
uniqueness of the solution of (\ref{fix}) on $[-T,0]$ for any
$T>0$.
Thus, we have a unique solution
of (\ref{fix}) on the whole real axis, i.e., $\mathbf{u}\in C(0,\R; L^q(I))$.
The integrand in (\ref{df-K}) is continuous as a map
$L^q(I)\to L^q(I)$. Thus, \eqref{df-K} and \eqref{fix} imply that $\mathbf{u}$ is
continuously differentiable and
we obtain a classical solution
of the IVP for (\ref{cKM+}) on the whole real axis. Finally,
since $K:\mathcal{M}\to\mathcal{M}$ is a uniform contraction (cf.~\eqref{contraction}),
which depends on $g$ continuously (cf.\eqref{df-K}), the fixed point is a continuous
function of $g$ as well (cf.~\cite[\S 1.2.6, Exercise~3]{Henry-Geom}).\\
$\qed$

\subsection{A priori estimates}\lbl{sec.apriori}

\begin{thm}\lbl{thm.apriori}
Let $\mathbf{u}(t)$ denote the solution of the IVP for
\eqref{cKM+} with $U\in L^1(U)$ and initial condition
$\mathbf{u}(0)=g\in L^\infty(I)$.
Then $\mathbf{u}\in C(\R; L^\infty(I))$ and
for any $T>0$, there exists $C>0$ depending on $T$ 
but not on $U$ such that
\be\lbl{Linfty-L2}
\|\mathbf{u} \|_{C(0,T; L^\infty(I))} \le C \|\mathbf{u}(0)\|_{L^\infty (I)}.
\ee
\end{thm}

\pf\; 
If $U\in L^1(I^2)$ and $\mathbf{u}(0)=g\in L^\infty(I)$ then the
contraction mapping argument used in the proof of
Theorem~\ref{thm.wellposed} yields exisitence of the unique solution 
$\mathbf{u}\in C^1(\R; L^\infty(I))$. Indeed, let 
$$
\mathcal{M}:=C(0,\tau; L^\infty(I)),\quad \mbox{for}\quad \tau:=(6L)^{-1}
$$
and consider the operator $K$ defined by \eqref{df-K}.
As before, we show that $K$ is a well defined contraction on $\mathcal{M}$. 

Indeed, for
any $\mathbf{u}, \mathbf{v}\in \mathcal{M},$ we have
\begin{multline*}
\left\| K\mathbf{u} - K\mathbf{v} \right\|_{\mathcal{M}} =
\max_{t\in [0,\tau]} \left\| K\mathbf{u}(t) - K\mathbf{v}(t)
\right\|_{L^\infty (I)} \\
\le \max_{t\in [0,\tau]} 
\left\| \int_0^t \left(\int_I U(\cdot,y)  
\left|  D\left( u(y,s)-u(\cdot,s)\right)- D\left( v(y,s)-v(\cdot,s)\right)
\right| dy   +|f(u(\cdot,s))-f(v(\cdot,s))|\right)  
ds \right\|_{L^\infty (I)}.
\end{multline*}
Using the Lipschitz continuity of $D$ and $f$ and the triangle 
inequality, we further obtain
\begin{multline*}
\left\| K\mathbf{u} - K\mathbf{v} \right\|_{\mathcal{M}} \le
\max_{t\in [0,\tau]} L\left\| \int_0^t \left(\int_I U(\cdot,y)
 \left|  u(y,s)-u(\cdot,s)-
v(y,s)+v(\cdot,s)
\right| dy   +|u(\cdot,s)-(v(\cdot,s)|\right)  ds \right\|_{L^\infty (I)}\\
\le \max_{t\in [0,\tau]} \left( 2L \int_I U(\cdot,y) +
L\right) \int_0^t\|u(\cdot,s)-v(\cdot,s)\|_{L^\infty (I)} ds
\le 3 L \tau \left\| \mathbf{u} - \mathbf{v} \right\|_{\mathcal{M}}\\
\end{multline*}
Recalling,  the definition of $\tau$, we arrive at
$$
\left\| K\mathbf{u} - K\mathbf{v} \right\|_{\mathcal{M}}  
\le {1\over 2} \left\| \mathbf{u} - \mathbf{v} \right\|_{\mathcal{M}}.
$$ 
Following the lines of the proof of Theorem~\ref{thm.wellposed}, 
it is straightforward to
show that the fixed point of \eqref{fix} is the unique solution
of the IVP for \eqref{cKM+}, $\mathbf{u}\in C^1(\R, L^\infty(I)),$
which depends continuously on the initial data $g\in L^\infty(I)$.

Denote
$$
m(t):= \| u(\cdot, t)\|_{L^\infty(I)}.
$$
From \eqref{cKM+}, using Lipschitz continuity of $D$ and $f,$ we have
\begin{eqnarray*}
\left| u(x,t) \right|& =& 
|g(x)| + L\int_0^t \left( \int_I U(x, y) \left( |u(y,s)-v(y,s)| + |u(x,s)-v(x,s)|\right) 
dy + |u(x,s)-v(x,s)| \right) ds\\ 
&\le& m(0)+ L( 2  \int_I U(x, y) dy +1) \int_0^t m(s)ds.
\end{eqnarray*}
Thus,
$$
m(t)\le m(0) + 3L \int_0^t m(s)ds.
$$
Since $\mathbf{u}\in C(0,T; L^\infty(I)),$ by Gronwall's inequality 
(cf.~\cite[Appendix B]{EvaPDE}), for any $t\in [0,T]$
$$
m(t)\le m(0)\left(1+3Lte^{3Lt}\right) \le C m(0), \quad C:=1+3LTe^{3LT},
$$
and \eqref{Linfty-L2} follows.\\
$\qed$

We will also use the following observation.

\begin{lem}\lbl{lem.dissipate}
Let $W\in L^2(I^2)$ be  a symmetric function and $u \in L^\infty(I)$.
Then
\be\lbl{dissipate+}
\int_{I^2} W(x,y) \left(u(y)-u(x)\right) u(x)dxdy={-1\over 2} \int_{I^2}
W(x,y) \left(u(y)-u(x)\right)^2 dxdy.
\ee
\end{lem}
\pf\;
Rewrite the left-hand side of \eqref{dissipate+} as
\begin{eqnarray}
\nonumber
\int_{I^2} W(x,y) \left(u(y)-u(x)\right) u(x)dxdy&=&
-\int_{I^2} W(x,y) \left(u(y)-u(x)\right)\left(u(x)-u(y)\right) dxdy\\
\lbl{square+}
&+&\int_{I^2} W(x,y) \left(u(y)-u(x)\right) u(y)dxdy.
\end{eqnarray}
Using the symmetry of $W(x,y)$, for the second term on the right-hand
side of \eqref{square+} we have
\be\lbl{switch-index+}
\int_{I^2} W(x,y) \left(u(y)-u(x)\right) u(y)dxdy=-\int_{I^2} W(x,y)
\left(u(y)-u(x)\right) u(x)dxdy.
\ee
After plugging \eqref{switch-index+} into \eqref{square+}, we obtain
\eqref{dissipate+}.\\
$\qed$

Next, we formulate the discrete counterparts of
Theorem~\ref{thm.apriori} and Lemma~\ref{lem.dissipate}.
To this end, consider an IVP for the semilinear discrete heat equation
\be\lbl{discrete-heat}
\dot u_{ni} = {1\over n} \sum_{j=1}^n V_{nij} D(u_{nj}-u_{ni}) + f(u_{ni}),\; i\in [n],
\ee
where $(V_{nij})$ is a nonnegative matrix with entries derived from
the graphon $W$ (see \eqref{def-Vnij}).

Let $u_n(t)=(u_{n1}(t), u_{n2}(t),\dots, u_{nn}(t))$ be a 
solution of \eqref{discrete-heat}. Denote
\be\lbl{discrete-norms}
\|u_n\|_{2,n}=\sqrt{n^{-1}\sum_{i=1}^n u_{ni}^2}\;\;\mbox{and}\;\;
\| u_n\|_{\infty,n} =\max_{i\in [n]} |u_{ni}|.
\ee

Recall that the discrete problem \eqref{discrete-heat} can be
rewritten as the nonlocal equation \eqref{average-step}.
By applying Theorem~\ref{thm.apriori} to \eqref{average-step},
we obtain the following theorem.
\begin{thm}\lbl{thm.discrete-apriori}
For the solution of the IVP for \eqref{discrete-heat}, we have
\be\lbl{discrete-Linfty-L2}
\max_{t\in [0,T]} \|u_n(t)\|_{\infty,n}  =C \|u_n(0)\|_{\infty, n} \; \forall n,
\ee
where $C>0$ depends on $T$ only.
\end{thm}

Finally, we state a discrete version of Lemma~\ref{lem.dissipate}. It can be derived from
Lemma~\ref{lem.dissipate}, or proved directly.

\begin{lem}\lbl{lem.discrete-dissipate}
Let $(W_{ij})$ be an $n\times n$ symmetric matrix. Then
for any $(\theta_1,\theta_2,\dots,\theta_n)\in\R^n$ 
\be\lbl{dissipate}
\sum_{i,j=1}^n W_{ij}(\theta_j-\theta_i) \theta_i={-1\over 2} \sum_{i,j=1}^n
W_{ij}(\theta_j-\theta_i)^2.
\ee
\end{lem}

\section{Averaging} \lbl{sec.averaging}
\setcounter{equation}{0}
In this section, we show that for large $n$ the solutions of the heat equation
\eqref{heat} on $\Gamma_n$ can be approximated by the solutions of the averaged
equation \eqref{aKM}, \eqref{average}.

For convenience, we rewrite the original and the averaged models. For the former
model, we plug in the expression for the mean degree $d_{ni}$ \eqref{exp-deg}
into \eqref{heat} to obtain
\be\lbl{new-heat}
\dot u_{ni}=n^{-1}\sum_{j=1}^n \eta_{nij}(u_{nj}-u_{ni}) + f(u_{ni}), \quad i\in [n],
\ee
where
\be\lbl{def-eta}
\eta_{nij}= \xi_{nij} \left(\rho_n  n^{-1}\sum_{j=1}^n \bar W_{nij}\right)^{-1}.
\ee
Recall the averaged model \eqref{aKM}:
\be\lbl{new-ave}
\dot v_{ni}=n^{-1}\sum_{j=1}^n  V_{nij} (v_{nj}-v_{ni}) +f(v_{ni}), \quad i\in [n],
\ee
where 
\be\lbl{def-Gni}
V_{nij}=\bar G_{ni}^{-1}  W_{nij},\quad G_{ni}:= n^{-1}\sum_{j=1}^n \bar W_{nij}.
\ee
Note that for fixed $i\in [n]$, $\{ \eta_{nij},\; j\in [n]\}$ are independent
random variables and
\be\lbl{Eeta}
\E \eta_{nij}= U_{nij},\quad i,j\in [n].
\ee
Below, we use the following weighted norm in $\R^n:$
\be\lbl{G-norm}
\|\psi_n\|_{G_n}:=\sqrt{ n^{-1} \sum_{i=1}^n G_{ni} \psi_{ni}^2}.
\ee
Here, we implicitly assume that   $n$ is large enough, so that
$\min_{i\in [n]} G_{ni}>0$ (cf.~\textbf{W-3)}).

We now formulate the main result of this section.
\begin{thm}\lbl{thm.ave}
Let $u_n(t)$ and $v_n(t)$ denote solutions of the IVP for
\eqref{new-heat} and \eqref{new-ave} respectively. Suppose that
the initial data for these problems satisfy 
\be\lbl{bdd-ic} 
\max\{ |u_n(0)|, |v_n(0)|\}\le C_1 \;\;\mbox{uniformly in}\;n \;
\mbox{and} \ee \be\lbl{close-ic} \lim_{n\to\infty}\left\|
v_n(0)-u_n(0)\right\|_{G_n}=0. \ee 
Then \be\lbl{ave}
\lim_{n\to\infty} \max_{t\in [0,T]}\left\|
v_n(t)-u_n(t)\right\|_{G_n} =0 \quad \mbox{a.s.}.
\ee
\end{thm}

For the proof of Theorem~\ref{thm.ave}, we will need the following lemma.
\begin{lem}\lbl{lem.Markov}
Let $T>0$ and $(a_{nij}(t))$ be an $n\times  n$ matrix, whose
entries depend on $t\in [0,T].$ Suppose \be\lbl{bdd-a} \sup_{t\in
[0,T]} \max_{(i,j)\in [n]^2}|a_{nij}(t)|\le C_2\;\; \forall n. \ee
Define  $Z_n(t)=(Z_{n1}(t), Z_{n2}(t), \dots, Z_{nn}(t))$, where
\be\lbl{def-Z} Z_{ni}(t)=n^{-1} \sum_{j=1}^n a_{nij} (t)
(\eta_{nij}-V_{nij}), \; i\in [n],\; t\in [0,T], \ee $\eta_{nij}$
are defined in \eqref{def-eta} (see \eqref{xi} for the definition
of $\xi_{nij}$.)

Then
\be\lbl{Z-as}
\lim_{n\to\infty} \sup_{t\in [0,T]} \|Z_n(t)\|_{G_n} =0\quad \mbox{a.s.}.
\ee
\end{lem}
\pf
Using the definitions \eqref{def-eta}, \eqref{def-Z}  and the bound in
\eqref{bdd-a}, for arbitrary $t\in [0,T],$ we have
\begin{eqnarray}
\nonumber
\E \sup_{t\in [0,T]} Z_{ni}(t)^2 &=& n^{-2} \E \left( \sup_{t\in [0,T]} \sum_{k,j=1}^n a_{nij}(t)
 a_{nik}(t) \left(\eta_{nij}-V_{nij}\right)
(\eta_{nik}-U_{nik})\right)\\
\nonumber &=&
n^{-2} \E  \sum_{j=1}^n \sup_{t\in [0,T]} a_{nij}(t)^2 \E (\eta_{nij}-V_{nij})^2=n^{-2}
\E\sum_{j=1}^n  \sup_{t\in [0,T]} a_{nij}(t) (\E\eta_{nij}^2-V_{nij}^2)\\
\lbl{2nd-m} &\le & \rho_n^{-1} n^{-2} \sum_{j=1}^n \sup_{t\in
[0,T]} a_{nij}(t)^2 V_{nij} \le C_2^2 \rho_n^{-1} n^{-2}
\sum_{j=1}^n V_{nij}= C_2^2 \rho_n^{-1} n^{-2}.
\end{eqnarray}
Let $\epsilon>0$ and $t\in [0,T]$ be arbitrary but fixed and denote the event
$$
A_n=\left\{ \sup_{t\in [0,T]} \| Z_n(t)\|_{G_n}  \ge \epsilon \right\}.
$$
By Markov's inequality, for arbitrary $\epsilon>0$
\be\lbl{Markov} 
\P\left(\sup_{t\in [0,T]} n^{-1}\sum_{i=1}^n G_{ni} Z_{ni}(t)^2 \ge
\epsilon\right)\le 
(\epsilon\rho_n n^2)^{-1} C_2^2 \left( n^{-2} \sum_{i,j=1}^n \bar
  W_{nij}\right). 
\ee 
where we used the definition of $G_{ni}$ in \eqref{def-eta}.
Since $\rho_nn\to \infty$ as $n\to\infty$ and 
$$
\lim_{n\to\infty} n^{-2} \sum_{i,j=1}^n \bar W_{nij} =\int_{I^2}
W(x,y)dxdy \quad 
(\mbox{cf.}~\eqref{W-2}),
$$
we have
$$
\sum_{n=1}^\infty \P(A_n)<\infty.
$$
By Borel--Cantelli Lemma \cite{Shir-Prob}, 
$\P(A_n \;\mbox{holds infinitely often}\;)=0$, i.e., 
$\sup_{t\in [0,T]} \|Z_{n}(t)\|_{G_n} \to 0$  a.s.  as $n\to\infty$.
This proves the lemma.\\
$\qed$

\pf (Theorem~\ref{thm.ave}) Denote $\phi_{ni}=u_{ni}-v_{ni}$.
By subtracting (\ref{aKM}) from (\ref{heat}), multiplying the result by
$n^{-1} G_{ni}\phi_{ni},$ and summing over $i\in [n]$, we obtain
\begin{eqnarray}\nonumber
2^{-1} {d\over dt} \|\phi_n\|^2_{G_n} &=& n^{-2}\sum_{i,j=1}^n \bar W_{nij}(\phi_{nj}-\phi_{ni}) \phi_{ni}
+ n^{-2} \sum_{i,j=1}^n G_{ni} (\eta_{nij}-V_{nij})(u_{nj}-u_{ni}) \phi_{ni}\\
\lbl{subtract}
&+& n^{-1}\sum_{i=1}^n G_{ni} \left[ f(u_{ni})-f(v_{ni})\right]\phi_{ni}.
\end{eqnarray}
By
Lemma~\ref{lem.discrete-dissipate}, the first term on the right--hand side of
\eqref{subtract} is nonpositive
\be\lbl{use-lemma}
\sum_{i,j=1}^n \bar W_{nij}(\phi_{nj}-\phi_{ni}) \phi_{ni} =-2^{-1}
\sum_{j=1}^n \bar W_{nij}(\phi_{nj}-\phi_{ni})^2\le 0.
\ee
Thus, using \eqref{use-lemma} and \eqref{Lip-f}, from \eqref{subtract}
we have
\be\lbl{pre-Gron}
2^{-1} {d\over dt} \|\phi_n\|^2_{G_n} \le
n^{-2} \sum_{i,j=1}^n G_{ni} (\eta_{nij}-V_{nij})(u_{nj}-u_{ni}) \phi_{ni}+L_f\|\phi_n\|_{G_n}^2.
\ee
Further, denote
$$
a_{nij}(t):=u_{nj}(t)-u_{ni}(t),\; Z_{ni}(t):=n^{-1}\sum_{j=1}^n
a_{nij}(t) (\eta_{nij}-V_{nij}), \; (i,j)\in [n]^2,\; t\in [0,T],
$$
and use $ab\le 2^{-1}(a^2+b^2)$ to obtain
\be\lbl{triangle+}
\left| n^{-2} \sum_{i,j=1}^n G_{ni} (\eta_{nij}-V_{nij})(u_{nj}-u_{ni}) \phi_{ni}\right|
\le 2^{-1} ( \|Z_n\|_{G_n}^2 +\|\phi_n\|_{G_n}^2).
\ee
Using \eqref{triangle+}, from \eqref{pre-Gron} we obtain
\be\lbl{pre-Gron-again}
{d\over dt} \|\phi_n\|^2_{G_n} \le  (2L_f+1)\|\phi_n\|_{G_n}^2 +\|Z_n\|_{G_n}^2.
\ee

Using \eqref{bdd-ic}, from Theorem~\ref{thm.discrete-apriori}, we have
$$
\max_{t\in [0,T]} \max_{(i,j)\in [n]^2}|a_{nij}(t)|\le C_3, \;
\forall n.
$$
By Lemma~\ref{lem.Markov}, with probability $1$,
for a given $\epsilon >0$
$$
\exists N_1(\epsilon):\quad \sup_{t\in [0,T]}\|Z_n\|_{G_n}^2\ge \epsilon^2/2,
$$
whenever $n\ge N_1$. Thus, for such $n,$
\be\lbl{pre-Gron-1}
 {d\over dt} \|\phi_n\|^2_{G_n} \le (2L_f+1) \|\phi_n\|_{G_n}^2 +\epsilon^2/2.
\ee
Thus, by Gronwall's inequality, we obtain
$$
\sup_{t\in[0,T]}\|\phi_n(t)\|^2_{G_n}\le \|\phi_n(0)\|^2_{G_n} e^{(2L+1)T} +{\epsilon^2\over 2(2L_f+1)}.
$$
Furthermore, by \eqref{close-ic}
$$
\exists N_2(\epsilon, T):\; \forall n\ge N_2 \quad \|\phi_n(0)\|^2_{G_n} e^{(2L_f+1)T}\le \epsilon^2/2.
$$
Thus, for $n\ge N(T,\epsilon):=N_1\vee N_2$ we have
$$
\sup_{t\in [0,T]}\|\phi_n(t)\|_{G_n}\le \epsilon.
$$
$\qed$

\section{The continuum limit}\lbl{sec.climit}
\setcounter{equation}{0}

Having justified averaging in \eqref{aKM}, our next goal is to show that the IVP for the averaged
equation \eqref{aKM} can be approximated by that for the
continuum limit \eqref{cKM}, \eqref{kernel}, subject to the initial condition
\be \lbl{clim-ic}
u(x,0) = g(x),
\ee
where $g\in L^\infty(I)$. To compare the solutions of the discrete
problem \eqref{aKM} and continuous equation \eqref{cKM} we
supply the former problem with the initial condition that is
consistent with \eqref{clim-ic}:
\be\lbl{aKM-ic}
v_{ni}(0)=n\int_{I_{ni}} g(x)dx, \; i\in [n].
\ee

Below, we construct a
finite-dimensional Galerkin approximation of \eqref{cKM} and \eqref{clim-ic}
and prove its convergence. In the next section, we compare solutions obtained by the
Galerkin's scheme with the solutions of the IVP for \eqref{aKM}.

Throughout this section, we assume that conditions
\textbf{W-1)}-\textbf{W-3)} and \eqref{Lip-f} hold.

\subsection{The Galerkin problem} \lbl{sec.galerkin}

Let $X=L^2(I)$,  define $K:X\to X$ by
\be\lbl{def-K}
[K(u)](x)=\int_I U(x,y) (u(y)-u(x))dy.
\ee
and rewrite \eqref{cKM} as follows
\begin{eqnarray}\lbl{Xode}
 \mathbf{u}^\prime &=& K(\mathbf{u}) + f(\mathbf{u}),\\
\lbl{Xode-ic}
\mathbf{u}(0) &=& g.
\end{eqnarray}
Recall that
$\mathbf{u}:\R\to X$ stands for the vector-valued function defined by $[\mathbf{u}(t)](x)=u(x,t)$ for
each $t\in\R$.

\begin{df}\lbl{weak-solution}
Function $\mathbf{u}\in H^1(0,T; X)$ is called a weak solution of the IVP
(\ref{Xode}), (\ref{Xode-ic})
on $[0,T]$ if
\be\lbl{weak-sol}
\left( \mathbf{u^\prime}(t)- K(\mathbf{u}(t)) -f(\mathbf{u}(t)), \mathbf{v}\right)=0\quad
\forall \mathbf{v}\in X
\ee
almost everywhere (a.e.) on $[0,T]$
and $\mathbf{u}(0)=g$.
\end{df}

To construct a finite-dimensional problem approximating
\eqref{weak-sol}, we introduce
$X_n=\operatorname{span}\{\phi_{ni}:\, i\in [n]\},$
a linear subspace of $X$. Here,
\be\lbl{def-phi}
\phi_{ni}(x)=\one_{I_{ni}}(x)=\left\{ \begin{array}{ll}
1, & x\in I_{ni},\\
0, & x \not\in I_{ni},
\end{array}\right.\; i\in [n].
\ee

Next, we construct the Galerkin approximation of the solution of
\eqref{Xode}, \eqref{Xode-ic}. To this end, we fix $n\in \N$ and
look for the approximate solution in the form 
\be\lbl{Gsol}
\mathbf{u}_n(t)=\sum_{i=1}^n u_{ni}(t)\mathbf{\phi}_{ni}. 
\ee 
The differentiable coefficients $u_{ni}(t),\; i\in [n],$ are
determined by projecting the original equation and the initial
condition on $X_n$: 
\be\lbl{project} 
\left(\mathbf{u}_n^\prime(t)- K(\mathbf{u}_n(t))-f(\mathbf{u}_n(t)),
\phi \right)=0\quad \forall \mathbf{\phi}\in X_n, 
\ee
\be\lbl{project-ic} \mathbf{u}_n(0)=P_{X_n} g=\sum_{i=1}^n
{(g,\phi_{ni})\over (\phi_{ni},\phi_{ni})} \phi_{ni}, \ee 
where
$P_{X_n}: X\to X_n$ stands for  the orthogonal projector  onto
$X_n.$ After plugging \eqref{Gsol} into \eqref{project} and
setting $\mathbf{v}=\phi_{ni},\; i\in [n]$, we arrive at the
following IVP  for the unknown coefficients $u_{ni}(t), \; i\in
[n]$:
\begin{eqnarray}\lbl{ode}
\dot u_{ni}(t)&=&n^{-1}\sum_{j=1}^n U_{nij} \left(
  u_{nj}(t)-u_{ni}(t)\right) + f(u_{ni}), \\
\lbl{ode-ic}
u_{ni}(0) &=& {(g,\phi_{ni})\over (\phi_{ni},\phi_{ni})} =n\int_{I_{ni}} g(x) dx.
\end{eqnarray}
Here,
\be\lbl{def-Unij}
U_{nij}  = n^2 \int_{I_{ni}\times I_{nj}} U(x,y) dxdy=
n^2\int_{I_{ni}\times I_{nj}} {W(x,y) \over \int_I W(x,z)dz} dxdy\le n.
\ee
Note that the right--hand side of \eqref{ode} is uniformly
Lipschitz continuous, which guarantees the existence of a unique
solution of the IVP \eqref{ode}, \eqref{ode-ic} on $\R$.

It will be convenient to have the Galerkin equation (\ref{ode})  rewritten as the integral equation
\be\lbl{finite-int}
{\p \over \p t} u_n(x,t) =\int_{I}  U_n(x,y) \left(u_n(y,t) -u_n(x,t)\right)dy + f(u_n(x,t)),
\ee
where $U_n$  and $u_n$ are  step functions
\begin{eqnarray}\lbl{def-Un}
U_n (x,y) &=& \sum_{i,j=1}^nU_{nij}\one_{I_{ni}\times I_{nj}}(x,y),\\
\nonumber
u_n(x,t)&=& \sum_{i=1}^n u_{ni}(t) \one_{I_{ni}}(x).
\end{eqnarray}

\subsection{Convergence of the Galerkin scheme} \lbl{sec.converge}

In this section, we show that the solutions of the Galerkin problems \eqref{project}, \eqref{project-ic},
$\mathbf{u}_n$, converge to $\mathbf{u}$, a unique weak solution of
\eqref{Xode}, \eqref{Xode-ic},  in the $L^2(0,T; X)$ norm as $n\to\infty$.

\begin{thm}\lbl{thm.converge} For any $T>0$,
there is a unique weak solution of \eqref{Xode}, \eqref{Xode-ic},
$\mathbf{u}\in H^1(0,T;X).$
The solutions of the Galerkin problems \eqref{project},
\eqref{project-ic}, $\mathbf{u}_n$ converge to
$\mathbf{u}$ in the $L^2(0,T; X)$ norm as $n\to\infty$.
\end{thm}
\pf
\begin{enumerate}
\item We shall first establish the following bounds for the
solutions $\mathbf{u}_n$ of the Galerkin problem \eqref{project},
\eqref{project-ic} that hold uniformly in $n$ 
\be\lbl{Linfty-X-X*}
\exists C_4=C_4(T, \|\mathbf{u}(0)\|_{L^\infty(I)}):\quad \max\{
\|\mathbf{u}_n \|_{C(0,T; L^\infty(I))}, \|\mathbf{u}_n\|_{C(0,T;
X)}, \|\mathbf{u}_n^\prime  \|_{C(0,T; X)} \}\le C_4. \ee 
The $L^\infty$-bound and, therefore, the $X$--bound follow from
Theorem~\ref{thm.apriori}. These bounds are uniform in $n$, 
because
$$
\|\mathbf{u}_n(0)\|_{L^\infty(I)}= \|P_{X_n} g\|_{L^\infty(I)}\le \|g\|_{L^\infty(I)}.
$$

To bound $\|\mathbf{u}^\prime_\mathbf{n}\|_{C(0,T;X)}$ we proceed as follows
$$
\left|(\mathbf{u_n^\prime}(t), \mathbf{v})\right| \le
\int_I \bar U_n(x,y) |u_n(x,t)-u_n(y,t)| |v(x,t)|dx dy+ \int_I |f(u_n(x))|| v(x)| dx.
$$
Using the $L^\infty$-bound for $\mathbf{u}_n$ \eqref{Linfty-X-X*}, the
continuity of $f$,  $\|\bar{U}_n\|_{L^2(I^2)}\le \|U\|_{L^2(I^2)} ,$ and the
triangle and Cauchy-Schwarz
inequalities, we  obtain
$$
\left|(\mathbf{u_n^\prime}(t), \mathbf{v})\right| \le
C_5(\|U\|_{L^2(I^2)} +C_6) \|\mathbf{v}\| \;\forall \mathbf{v}\in X.
$$
Thus, \be\lbl{H1pt} \|\mathbf{u_n^\prime}(t)\|\le C_6, t\ge 0. \ee
uniformly in $n$. \item Estimates in \eqref{Linfty-X-X*} imply
\begin{eqnarray}
\lbl{L20T}
\|\mathbf{u}_n\|_{L^2(0,T;X)} &\le& C_4,\\
\lbl{u-Lip} \| \mathbf{u}_n(t+h)-\mathbf{u}_n(t)\|_X& \le& C_4
|h|,
\end{eqnarray}
respectively. From \eqref{u-Lip}, we further have \be\lbl{L2cont}
\int_0^T\| \mathbf{u}_n(t+h)-\mathbf{u}_n(t)\|_X^2dt\le C_4^2T
h^2. \ee From \eqref{L20T}  and (\ref{L2cont}), using the
Frechet--Kolmogorov theorem (cf. \cite{Yosida-Analysis}), we see
that $(\mathbf{u}_n)$ is precompact in $L^2(0,T;X)$. Thus, one can
select a subsequence $(\mathbf{u}_{n_k})$ that converges to
$\mathbf{u}\in L^2(0,T;X)$. \item Likewise, integrating both sides
of \eqref{H1pt} from $0$ to $T$, we obtain
$$
\|\mathbf{u_n^\prime}\|_{L^2(0,T;X)}\le C_6\sqrt{T}
$$
uniformly in $n$. Thus, $(\mathbf{u}^\prime_{n_k})$ is weakly
precompact in $L^2(0,T;X)$, and one can select  a subsequence
$(\mathbf{u}^\prime_{n_{k^\prime}})$ that weakly converges to
$\mathbf{w}\in L^2(0,T;X)$ and strongly converges to
$\mathbf{u}^\prime\in L^2(0,T;X)$. Clearly,
$\mathbf{w}=\mathbf{u^\prime}$. Indeed, taking $\phi \in C^1(0,T;
X)$ with compact support in $(0,T)$ and using integration by
parts, we obtain \be\lbl{byparts} \int_0^T
\mathbf{u}_{n_{k^\prime}}^\prime(t) \phi (t) dt = -\int_0^T
\mathbf{u}_{n_{k^\prime}}(t) \phi^\prime (t)  dt. \ee By sending
$k^\prime \to\infty$ in (\ref{byparts}), we see that
$\mathbf{u}^\prime\in H^1(0,T; X)$ and
$\mathbf{u}^\prime=\mathbf{w}.$ \item Next, we show that
$\mathbf{u}$ is a unique weak
  solution of \eqref{weak-sol}
satisfying $\mathbf{u}(0)=\mathbf{g}$.
This follows from a standard argument (see, e.g., \cite[Theorem~7.1.3]{EvaPDE}).

Fix $N\in\N$ and choose a function of the form
\be\lbl{separable}
\mathbf{v}(t)=\sum_{j=1}^{N} d_j(t) \mathbf{\phi_{nj}},
\ee
where $d_j(t)$ are continuously differentiable functions and $\mathbf{\phi}_{nj}$ are
defined in \eqref{def-phi}.
Next, we multiply \eqref{project} with
$n>N$ and $\mathbf{\phi}:=\phi_{nj}$ by $d_j(t),$ sum over $j$, and integrate the result
from $0$ to $T$ to obtain
$$
\int_0^T
(\mathbf{u^\prime}_n(t)-K(\mathbf{u}_n(t))-f(\mathbf{u}_n(t)),
\mathbf{v}(t))dt=0.
$$
Passing to the limit along $n=n_k$, we have \be\lbl{w-limit}
\int_0^T (\mathbf{u}^\prime(t)-K(\mathbf{u}(t))
-f(\mathbf{u}(t)),\mathbf{v}(t)) dt=0. \ee This equality holds for
an arbitrary $\mathbf{v}$ as in \eqref{separable}. Since functions
of this form for $N\in\N$ are dense in $L^2(0,T;X)$, we conclude
that \eqref{w-limit} holds for all $\mathbf{v}\in L^2(0,T;X)$.
Therefore, \be\lbl{weak-equality} (\mathbf{u^\prime}-K(\mathbf{u})
-f(\mathbf{u}),\mathbf{v})=0 \quad \forall \mathbf{v}\in
L^2(0,T;X) \ee a.e. on $[0,T]$.
\item To show that  $\mathbf{u}$ is a weak solution of
\eqref{cKM}, \eqref{clim-ic}, it remains to verify
$\mathbf{u}(0)=\mathbf{g}$. To this end, we choose $\mathbf{v}\in
C^1(0,T;X)$ vanishing at $t=T$ as a test function in
(\ref{weak-sol}) and integrate by parts to obtain \be\lbl{ic1}
-\int_0^T\left(\mathbf{u}(t), \mathbf{v^\prime}(t)\right)dt=
\int_0^T \left( K(\mathbf{u}(t)) + f(\mathbf{u}(t)),
\mathbf{v}(t)\right)dt+\left( \mathbf{u}(0), \mathbf{v}(0)\right).
\ee Using the same test functions in (\ref{project}), we have
\be\lbl{ic2} -\int_0^T\left(\mathbf{u}_{n_k}(t),
\mathbf{v^\prime}(t)\right)dt= \left(K(\mathbf{u}_{n_k}(t)) +
f(\mathbf{u}_{n_k}(t)), \mathbf{v}(t)\right)dt+ \left(
\mathbf{u}_{n_k}(0), \mathbf{v}(0)\right). \ee Passing to the
limit in (\ref{ic2}) yields \be\lbl{ic2-limit}
-\int_0^T\left(\mathbf{u}(t), \mathbf{v^\prime}(t)\right)dt=
\int_0^T \left(K(\mathbf{u}(t)) + f(\mathbf{u}(t)),
\mathbf{v}(t)\right)dt + \left( \mathbf{g}, \mathbf{v}(0)\right).
\ee Comparing the limiting equation (\ref{ic2-limit}) with
(\ref{ic1}) we conclude that $\mathbf{u}(0)=\mathbf{g}$. Thus,
$\mathbf{u}$ is a weak solution of \eqref{Xode}. 
\item To show
that the just constructed weak solution is unique, suppose that
there is another solution 
\be\lbl{another} \mathbf{w^\prime}=
K(\mathbf{w})+f(\mathbf{w}) 
\ee 
satisfying the same initial
condition $\mathbf{w}(0)=\mathbf{g}$. 
Denote $\mathbf{\xi}=\mathbf{u}-\mathbf{w}.$
By subtracting the
(\ref{another}) from (\ref{Xode}), multiplying both sides by
$G(x)\mathbf{\xi}$ and integrating over
$I,$  we have
$$
\frac{1}{2}{d\over dt} \| \sqrt{G}\xi(\cdot,t)\|_X^2=\int_{I^2} W(x,y)
\left(\xi(y,t)-\xi(x,t)\right)\xi(x,t)dydx
+ \int_I G(x) \left(f(u(x,t))-f(w(x,t)\right)\xi(x,t)dx
$$
Using Lemma~\ref{lem.dissipate} and Lipschitz continuity of $f$, we obtain
$$
{d\over dt} \| \sqrt{G}\xi(\cdot,t)\|_X^2\le L \|\sqrt{G} \xi(\cdot,t)\|_X^2.
$$
Since $G$ is strictly positive on $I$ (cf.~\eqref{int-below}),
from the last inequality and $\mathbf{\xi}(0)=0,$ we conclude that
$\mathbf{u}(t)=\mathbf{w}(t)$
for all $t\in[0,T]$. This proves uniqueness.
\item
The uniqueness of the weak solution entails $\mathbf{u}_n\to\mathbf{u}$
as $n\to\infty$. Indeed, suppose on the contrary that there exists a subsequence
$\mathbf{u}_{n_l}$, which is not converging to $\mathbf{u}$. Then for a given $\epsilon>0$
one can select a subsequence $\mathbf{u}_{n_{l_i}}$ such that
$$
\|\mathbf{u}_{n_{l_i}} -\mathbf{u}\|_{L^2(0,T;X)}>\epsilon \;\forall i\in\N.
$$
However, $(\mathbf{u}_{n_{l_i}})$ is precompact in $L^2(0,T,X)$ and contains a
subsequence converging to a weak solution of \eqref{Xode}, which must be
$\mathbf{u}$ by uniqueness. Contradiction.
\end{enumerate}
$\qed$

\subsection{Approximation}\lbl{sec.approximate}

It remains to estimate the difference between the solutions of the
averaged equation \eqref{aKM} and that of the Galerkin problem
\eqref{ode}.
The key is the estimate for the $L^4$--norm of the difference between
$V_n$ and $U_n$ (see \eqref{W-step} and \eqref{def-Un}),
the kernels used in the averaged and the Galerkin
problems respectively.

\begin{lem}\lbl{lem.martingale}
\be\lbl{Un-Vn}
\|U_n-V_n\|_{L^4(I^2)} \to 0,\quad\mbox{as}\; n\to\infty.
\ee
\end{lem}
\pf
First, we show that $U_n, n\in\N,$ form a sequence of $L^4$--bounded martingales 
\cite{Williams-Prob-Mart}. 
To this end, we consider a probability space
$(I^2, \mathcal{B}(I^2),\lambda)$ with $I^2$ as a sample space equipped with the
$\sigma$--algebra of Borel sets, and the Lebesgue measure as probability. Let
$\mathcal{A}_n$ denote the algebra of subsets of $I^2$ generated by the sets
$I_{ni}\times I_{nj}$, $(i,j)\in [n]^2.$ Then $U_n$ can be represented as the conditional
expectation
$$
U_n=\E (U|\mathcal{A}_n), \quad n\in\N.
$$
Since $U\in L^4(I^2)$ (cf.~\textbf{W-4)} and \eqref{int-below}),  
the $L^p$--Martingale Convergence Theorem yields
\be\lbl{Un-to-U}
U_n\to U \quad \mbox{a.e.
and in}\;  L^4(I^2)\;\mbox{as}\; n\to\infty.
\ee

Next, we turn to functions $V_n, n\in \N$ (cf.~ \eqref{W-step}):
\begin{eqnarray}
\nonumber
V_n(x,y)&=&\sum_{i,j=1}^n {\bar W_n(x_{ni}, x_{nj}) \over n^{-1} \sum_{k=1}^n \bar W_n(x_{ni}, x_{nk})}
\one_{I_{ni}\times I_{nj}} (x,y)\\
\nonumber
&=&
{\sum_{i,j=1}^n \bar W_n(x_{ni}, x_{nj}) \one_{I_{ni}\times I_{nj}}
  (x,y) \over \sum_{i=1}^n n^{-1} \sum_{k=1}^n \bar W_n(x_{ni},
  x_{nk}) \one_{I_{ni}} (x) }
\\
\lbl{rewtite-Vn}
&=:& {P_n(x,y)\over Q_n(x)}.
\end{eqnarray}

From \eqref{Wn} and \textbf{W-1)}, we have $P_n\to W$ a.e. on
$I^2$. Likewise,
by \eqref{W-3}, 
$$
Q_n= \int_I W(\cdot,z)dz(1+\delta_n)\ge \nu >0, \;\mbox{as}\; n\to\infty
$$
uniformly on any closed interval lying in $(0,1)$.
Thus,
$
{P_n\over Q_n}\to U
$
a.e. on $I^2$ as $n\to\infty.$ Furthermore, by \eqref{Wn} and \textbf{W-3)},
$$
\left| {P_n\over Q_n}\right|\le {W\over \nu}.
$$
Since $V_n, U\ge 0,$
$$
\left(V_n-U\right)^4\le 4(V_n^4+U^4)\le {4W^2\over \nu^4}.
$$
Thus, $V_n-U\to 0$ a.e. on $I^2$ and 
$$
\left| V_n-U\right|\le \sqrt{2}\nu^{-1} W\in L^4(I^2).
$$
By the Dominated Convergence Theorem $V_n\to U$ in $L^4(I^2)$.  From this and
\eqref{Un-to-U}, we conclude 
$$
\|U_n-V_n\|_{L^4(I^2)}\to 0, \quad n\to\infty.
$$
$\qed$

\begin{lem}\lbl{lem.approximate}
For any  $T>0$, solutions of the IVPs for \eqref{aKM} and \eqref{ode}  satisfy
\be\lbl{compare-ivps}
\lim_{n\to\infty}\max_{t\in [0,T]} \| u_n(t)-v_n(t)\|_{G_n} =0,
\ee
provided
\be\lbl{ic-converge}
\lim_{n\to\infty} \| u_n(0)-v_n(0)\|_{G_n} =0,
\ee
\end{lem}
\pf
Denote $\phi_{ni}:=u_{ni}-v_{ni}, \; i\in [n].$ By subtracting (\ref{aKM}) from (\ref{ode}), multiplying
the result by $n^{-1}G_{ni} \phi_{ni}$ (see \eqref{def-eta} for the definition of $G_{ni}$)
 and summing over $i\in [n]$, we obtain
\begin{eqnarray}\nonumber
2^{-1} {d\over dt} \|\phi_n\|^2_{G_n} &=& n^{-2}\sum_{i,j=1}^n \bar W_{nij}(\phi_{nj}-\phi_{ni}) \phi_{ni}
+ n^{-2} \sum_{i,j=1}^n G_{ni}(\bar U_{nij}-V_{nij})(u_{nj}-u_{ni}) \phi_{ni}\\
\lbl{as-before}
&+& n^{-1}\sum_{i=1}^nG_{ni} \left[ f(u_{ni})-f(v_{ni})\right]\phi_{ni}.
\end{eqnarray}
As before, we use Lemma~\ref{lem.discrete-dissipate} and \eqref{Lip-f} to obtain
\begin{eqnarray}\lbl{use-dissipate}
&&n^{-2}\sum_{i,j=1}^n \bar W_{nij}(\phi_{nj}-\phi_{ni}) \phi_{ni}\le 0,\\
\lbl{use-Lip-f}
&&\left|n^{-1}\sum_{i=1}^n G_{ni} \left[ f(u_{ni})-f(v_{ni})\right]\phi_{ni}\right|\le L_f\|\phi_n\|_{G_n}^2.
\end{eqnarray}
Using $\max_{t\in [0,T]}\|u_n(t)\|_{\infty,n}\le C_7$ (cf. Theorem~\ref{thm.discrete-apriori}),
we estimate
$$
\left|n^{-2} \sum_{i,j=1}^n G_{ni} (U_{nij}-V_{nij})(u_{nj}-u_{ni}) \phi_{ni}\right|\le
C_8(\Delta_n(W)+\|\phi_n\|_{G_n}^2),
$$
where
$$
\Delta_n(W):= {1\over n^2}\sum_{i,j=1}^n G_{ni} (U_{nij}-V_{nij})^2.
$$
Further, using the Cauchy-Schwarz inequality, we have
\be\lbl{Delta-W}
\Delta_n(W)\le \left(n^{-1}\sum_{i=1}^n G_{ni}^2 \right)^{1/2} \|U_n-V_n\|^2_{L^4(I^2)} 
\ee
Recalling the definition of $G_{ni}$ and using the Cauchy-Schwarz inequality again, we obtain
\be\lbl{use-W4}
n^{-1}\sum_{i=1}^n G_{ni}^2 = n^{-1} \sum_{i=1}^n 
\left( n^{-1} \sum_{j=1}^n \bar W_{nij} \right)^2\le n^{-2} \sum_{i,j=1}^n \bar W_{nij}^2.
\ee
Using \eqref{W-4}, \eqref{use-W4}, and Lemma~\ref{lem.approximate},
we obtain
\be\lbl{Darbu}
\lim_{n\to\infty}\Delta_n(W)=0.
\ee
The combination of \eqref{as-before}-\eqref{Delta-W} yields
$$
{d\over dt}\|\phi_n\|^2_{G_n} \le 2(C_8+L_f) \|\phi_n\|^2_{G_n}+2C_7 \Delta_n(W).
$$
By Gronwall's inequality,
\be\lbl{Gronwall-again}
\max_{t\in [0,T]} \|\phi_n(t)\|^2_{G_n} \le \left(\|\phi_n(0)\|^2_{2,n} +{C_8\over C_8+L_f} 
\Delta_n(W)\right) \exp\{(C_8+L_f)T\}.
\ee
The right hand side in \eqref{Gronwall-again} tends to $0$ as $n\to\infty$, as follows from
\eqref{ic-converge} and \eqref{Darbu}. This proves the lemma.\\
$\qed$

Theorem~\ref{thm.main} now follows from Theorems~\ref{thm.ave},
\ref{thm.converge} and Lemma~\ref{lem.approximate}.

\section{Discussion}\lbl{sec.discussion}
\setcounter{equation}{0}
The analysis in the preceding sections justifies the continuum limit \eqref{cKM} for the 
semilinear heat equation \eqref{heat} on sparse W--random graphs. In conclusion,
we outline several extensions of this work to certain nonlinear models, which are of 
interest in applications.

\subsection{The Kuramoto model}

The analysis in Sections~\ref{sec.averaging} and \ref{sec.climit} can
be extended to cover the following nonlinear heat equation:
\be\lbl{KM}
\dot u_{ni}= n^{-1}\sum_{ j: \{ i,j \} \in E(\Gamma_n)} D(u_{nj}-u_{ni}) + f(u_{ni}),
\ee
where $D$ and $f$ are Lipschitz continuous functions (cf.~\eqref{Lip}, \eqref{Lip-f}).
In addition, we assume that $D$ is an odd function satisfying the sign condition
\be\lbl{sign-D}
uD(u)\ge 0.
\ee
Both conditions hold for the original Kuramoto model with $D(u)=\sin u$. 

Under the above assumptions on $f$ and $D,$ we can justify the
continuum limit for \eqref{KM}.

\begin{thm}\lbl{thm.KM}
Let $g\in L^\infty(I)$ and $T>0$ be arbitrary. Denote the solutions of \eqref{KM} and
\eqref{nlin} subject to the initial conditions \eqref{ic-heat} and \eqref{ic-clim} by
$u_{ni}(t), \; i\in [n]$ and $u(x,t)$ respectively. 

Then with probability $1,$  
$$
\lim_{n\to 0}\sup_{t\in [0,T]} \|u_n(\cdot,t)- u(\cdot, t)\|_{L^2(I)} =0,
$$
where 
$$
u_n(x,t):=\sum_{i=1}^n u_{ni}(t) \one_{I_{ni}}(x).
$$
\end{thm}

For the proof of Theorem~\ref{thm.KM}, one needs the following modification of 
Lemma~\ref{lem.dissipate}.

\begin{lem}\lbl{lem.sym}
Let $W\in L^2(I)$ be a symmetric function and $D$ be an odd symmetric continuous 
function. Then for any $u\in L^\infty(I)$,
$$
\int_{I^2} W(x,y) D(u(y)-u(x)) u(x) dxdy= -2^{-1} \int_{I^2} W(x,y) D(u(y)-u(x)) (u(y)-u(x)) dxdy.
$$
If, in addition, $W\ge 0$ and $D$ satisfies \eqref{sign-D}, then
$$
\int_{I^2} W(x,y) D(u(y)-u(x)) u(x) dxdy\le 0.
$$
\end{lem}

With Lemma~\ref{lem.sym} in hand, the proofs of the statements
in Sections~\ref{sec.averaging} and \ref{sec.climit} can be translated to the nonlinear 
equation \eqref{KM} with minor changes.

\subsection{An alternative scaling and other graph models}
If the diffusion term is scaled by $n\rho_n$ instead of $d_{ni}=O(n\rho_n)$ as 
in \eqref{heat}, the formal derivation of the continuum limit yields
\be\lbl{alt-cont}
{\p\over \p t} u (x,t)=\int_I W(x,y) D\left( u(y,t)-u(x,t)\right)dy + f(u(x,t)).
\ee
Here, the kernel is $W$ instead of $U$ (cf.~\eqref{kernel}).
In particular, for the Kuramoto model on the power law family of graphs, the alternative scaling yields
\be\lbl{alt-KM}
{\p\over \p t} u (x,t)=x^{-\alpha}\int_I y^{-\alpha} \sin\left( u(y,t)-u(x,t)\right)dy.
\ee 
The presence of the $x$--dependent factor on the right--hand side of \eqref{alt-KM} has interesting
implications for the spatial patterns generated by the Kuramoto model. In particular,
it is responsible for the existence of the chimera-like patterns in the Kuramoto model
with repulsive coupling on power law graphs (cf.~\cite{MedTan16}).

The proof of existence of the strong solution of the IVP in Section~\ref{sec.ivp} does not
cover the equation \eqref{alt-cont}, because it relies on condition \eqref{norm}, which 
does not hold for $W$ in general (see \eqref{triangle}). However, one can show the existence
of the weak solution for the IVP for \eqref{alt-cont} (cf.~Definition~\ref{weak-solution})
by constructing it as the limit of solutions of the Galerkin problems following the lines
of the analysis in \S\ref{sec.converge}.

Likewise, there are many different ways how to define a convergent
family of sparse random graphs. Instead of \eqref{Pedge} one could, 
for example, define the probability for a given pair of nodes to
belong to the edge set using  averaging:
\be\lbl{alt-graph}
\P\left(\{i,j\}\in E(\Gamma_n)\right) =\rho_n n\int_{I_{ni}\times
  I_{nj}} W(x,y) dxdy.
\ee
The analysis of this paper can be used to justify the continuum limit for
coupled systems on $\{\Gamma_n\}$ defined by \eqref{alt-graph}.

\vskip 0.2cm
\noindent
{\bf Acknowledgements.}
This work was supported in part by the NSF DMS 1412066 (GM).

\bibliographystyle{amsplain}

\begin{thebibliography}{10}

\bibitem{BarAlb99}
A.-L. Barab{\'a}si and A.~Albert, \emph{Emergence of scaling in random
  networks}, Science \textbf{286} (1999), 509--512.

\bibitem{Biggs}
N.~Biggs, \emph{{A}lgebraic {G}raph {T}heory}, second edition ed., Cambridge
  University Press, 1993.

\bibitem{Bogachev-MT}
V.~I. Bogachev, \emph{Measure theory. {V}ol. {I}, {II}}, Springer-Verlag,
  Berlin, 2007. \MR{2267655 (2008g:28002)}

\bibitem{BCCZ-I-14}
C.~{Borgs}, J.~T. {Chayes}, H.~{Cohn}, and Y.~{Zhao}, \emph{{An $L^p$ theory of
  sparse graph convergence I: limits, sparse random graph models, and power law
  distributions}}, ArXiv e-prints (2014).

\bibitem{BCCZ-II-14}
\bysame, \emph{{An $L^p$ theory of sparse graph convergence II: {L}{D}
  convergence, quotients, and right convergence}}, ArXiv e-prints (2014).

\bibitem{CroHoh93}
M.C. Cross and P.C. Hohenberg, \emph{Pattern formation out of equilibrium},
  Rev. Mod. Phys. \textbf{65} (1993), 851--1112.

\bibitem{CucSma07}
F.~Cucker and S.~Smale, \emph{Emergent behavior in flocks}, IEEE Trans.
  Automat. Control \textbf{52} (2007), no.~5, 852--862. \MR{2324245
  (2008h:91132)}

\bibitem{DorBul12}
F.~Dorfler and F.~Bullo, \emph{Synchronization and transient stability in power
  networks and non-uniform {K}uramoto oscillators}, SICON \textbf{50} (2012),
  no.~3, 1616--1642.

\bibitem{EvaPDE}
L.C. Evans, \emph{Partial {D}ifferential {E}quations}, AMS, 2010.

\bibitem{Henry-Geom}
D.~Henry, \emph{Geometric theory of semilinear parabolic equations}, Lecture
  Notes in Mathematics, vol. 840, Springer-Verlag, Berlin-New York, 1981.
  \MR{610244 (83j:35084)}

\bibitem{Kur84}
Y.~Kuramoto, \emph{Cooperative dynamics of oscillator community}, Progress of
  Theor. Physics Supplement (1984), 223--240.

\bibitem{LovGraphLim12}
L.~Lov{\'a}sz, \emph{Large networks and graph limits}, AMS, Providence, RI,
  2012.

\bibitem{Med12}
G.S. Medvedev, \emph{Stochastic stability of continuous time consensus
  protocols}, SIAM Journal on Control and Optimization \textbf{50} (2012),
  no.~4, 1859--1885.

\bibitem{Med14a}
\bysame, \emph{The nonlinear heat equation on dense graphs and graph limits},
  SIAM J. Math. Anal. \textbf{46} (2014), no.~4, 2743--2766. \MR{3238494}

\bibitem{Med14b}
\bysame, \emph{The nonlinear heat equation on {W}-random graphs}, Arch. Ration.
  Mech. Anal. \textbf{212} (2014), no.~3, 781--803. \MR{3187677}

\bibitem{MedTan16}
G.S. Medvedev and X.~Tang, \emph{The {K}uramoto model on power law graphs}, in
  preparation.

\bibitem{MotTad14}
S.~Motsch and E.~Tadmor, \emph{Heterophilious dynamics enhances consensus},
  SIAM Rev. \textbf{56} (2014), no.~4, 577--621. \MR{3274797}

\bibitem{Shir-Prob}
A.N. Shiryaev, \emph{Probability}, Springer, 1996.

\bibitem{Williams-Prob-Mart}
D.~Williams, \emph{Probability with martingales}, Cambridge Mathematical
  Textbooks, Cambridge University Press, Cambridge, 1991. \MR{1155402
  (93d:60002)}

\bibitem{Yosida-Analysis}
K{\=o}saku Yosida, \emph{Functional analysis}, Classics in Mathematics,
  Springer-Verlag, Berlin, 1995, Reprint of the sixth (1980) edition.
  \MR{1336382 (96a:46001)}

\end{thebibliography}
\def\cprime{$'$}
\providecommand{\bysame}{\leavevmode\hbox to3em{\hrulefill}\thinspace}
\providecommand{\MR}{\relax\ifhmode\unskip\space\fi MR }
\providecommand{\MRhref}[2]{%
  \href{http://www.ams.org/mathscinet-getitem?mr=#1}{#2}
}
\providecommand{\href}[2]{#2}

\end{document}